\documentclass[10pt]{article}
\usepackage[tbtags]{amsmath}
\usepackage{epsfig,amstext,amssymb,amsthm,latexsym}
\usepackage{amssymb,color}

\newcommand{\nwc}{\newcommand}
\nwc{\COM}[1]{}

\definecolor{c20}{rgb}{0.,0.7,0.}
\definecolor{c30}{rgb}{0.,0.,1.}
\definecolor{c40}{rgb}{1,0.1,0.7}
\definecolor{c50}{rgb}{1,0,0}

\catcode`\@=11 \@addtoreset{equation}{section} \catcode`\@=12

\allowdisplaybreaks[4] 
\newtheorem{theo}{Theorem}[section]
\newtheorem{sat}[theo]{Proposition}
\newtheorem{de}[theo]{Definition}
\newtheorem{lem}[theo]{Lemma}

\newtheorem{korr}[theo]{Corollary}
\newtheorem{remark}[theo]{Remark}
\newtheorem{exxa}[theo]{Example}

\newcommand{\nelem}[1]{{Lemma \ref{#1}}}
\newcommand{\neprop}[1]{{Proposition \ref{#1}}}
\newcommand{\netheo}[1]{{Theorem \ref{#1}}}

\newcommand{\kb}[1]{\boldsymbol{#1}}
\newcommand{\vk}[1]{\kb{#1}}

\def\FRE{\mbox{Fr\'{e}chet }}

\def\X{\vk{X}}

\newcommand{\ve}{\varepsilon}

\newcommand{\abs}[1]{\lvert #1 \rvert}

\newcommand{\E}[1]{\mbox{\rm$\vk{E}$}\{#1\}}

\newcommand{\pk}[1]{\mbox{\rm$\vk{P}$} \{#1\} }

\newcommand{\R}{\!I\!\!R}

\newcommand{\N}{\!I\!\!N}
\newcommand{\inr}{\in \R}

\newcommand{\inn}{\in \N}

\newcommand{\ldot}{,\ldots,}

\newcommand{\limit}[1]{\lim_{#1 \to   \infty}}

\newcommand{\todis}{\stackrel{d}{\to}}

\newcommand{\equaldis}{\stackrel{d}{=}}
\newcommand{\eqd}{\stackrel{d}{=}}

\newcommand{\BQN}{\begin{eqnarray}}
\newcommand{\EQN}{\end{eqnarray}}
\newcommand{\BQNY}{\begin{eqnarray*}}
\newcommand{\EQNY}{\end{eqnarray*}}
\newcommand{\BS}{\begin{sat}}
\newcommand{\ES}{\end{sat}}
\newcommand{\BL}{\begin{lem}}
\newcommand{\EL}{\end{lem}}
\newcommand{\BT}{\begin{theo}}
\newcommand{\ET}{\end{theo}}
\newcommand{\BK}{\begin{korr}}
\newcommand{\EK}{\end{korr}}
\newcommand{\BD}{\begin{de}}
\newcommand{\ED}{\end{de}}
\newcommand{\BIT}{\begin{itemize}}
\newcommand{\EIT}{\end{itemize}}
\newcommand{\BDI}{\begin{description}}
\newcommand{\EDI}{\end{description}}

\newcommand{\BEX}{\begin{exxa}}
\newcommand{\EEX}{\end{exxa}}

\newcommand{\QED}{\hfill $\Box$}
\newcommand{\IF}{\infty}
\def\kal#1{{\cal{ #1}}}

\newcommand{\prooftheo}[1]{ \textsc{Proof of Theorem} \ref{#1} }
\newcommand{\proofprop}[1]{\textsc{Proof of Proposition} \ref{#1}}
\newcommand{\prooflem}[1]{\textsc{Proof of Lemma} \ref{#1}}

\def\SI{\Sigma}

\def\X{\vk{X}}

\def\HAB{H_{\alpha, \beta}}
\def\hab{h_{\alpha, \beta}}

\def\OPJ{ \kal{J}}

\def\Bab{B_{\alpha,\beta}}
\def\Fab{H_{\alpha,\beta}}

\def\gatrho{\alpha_{\widetilde \rho}}
\def\garho{\alpha_{\rho}}

\def\tRO{\widetilde \rho}

\begin{document}

\centerline{\Large On Beta-Product Convolutions}

        \vskip 1.8 cm
        \centerline{\large Enkelejd Hashorva \footnote{
Department of Actuarial Science, Faculty of Business and Economics, University of Lausanne, 
Extranef\\ UNIL-Dorigny 
1015 Lausanne, Switzerland, 
 email: enkelejd.hashorva@unil.ch}
}

      \centerline{\today{}}

 {\bf Abstract:} Let $R$ be a positive random variable independent of  $S$ which is beta distributed. In this paper we are interested  on the relation between $R$ and $RS$. For this model we derive first some distributional properties,
 and then investigate the lower tail asymptotics of $RS$ when $R$ is regularly varying at 0, and vice-versa.
 Our first application concerns  the asymptotic behaviour  of the componentwise sample minima related to an elliptical distributions. Further, we derive the lower tails asymptotic of the aggregated risk for bivariate polar distributions.

{\it AMS 2000 subject classification:} Primary 60F05; Secondary 60G70.\\

{\it Key words and phrases:} Weyl fractional-order
integral operator; Williamson $d$-transform; random scaling; elliptical distribution; polar distributions;
product convolution; asymptotics of sample minima; lower tail asymptotics; risk aggregation.

\section{Introduction}
Let $R$ and $S$ be two independent positive random variables. In this paper we consider the random scaling model
\BQN \label{eq:ab}
 W &\equaldis  & R S,
\EQN
with $W$ the scaled version of $R$ and $S\in (0,1)$ almost surely ($\equaldis$ stands for equality of the distribution functions).
 In order to derive distributional properties of $W$  we need to specify the distribution function of $S$; a tractable instance with various applications is the tractable case that $S$ is a beta distributed random variable.\\

In a financial or insurance framework, the random scaling model \eqref{eq:ab} appears naturally with  $W$ the  deflated risk arising from some loss or investment $R$ which is independent from the random scaling/deflating factor $S$.
Other prominent applications in the literature concern modeling of
network data (see e.g., D'Auria  and Resnick (2006, 2008));
 random difference equations (see e.g., Mikosch and Konstantinides (2004), Denisov  and  Zwart (2007));
 insurance and finance applications (see e.g., Tang and Tsitsiashvili (2003, 2004), Tang (2006, 2006),
 Piterbarg et al.\ (2009), Liu and Tang (2010), Tang and Vernic  (2010), Zhang (2010));
  approximation of multivariate distributions (see e.g., Hashorva (2007), Charpentier and Segers (2009), McNeil and Ne$\check{s}$lehov\'{a} (2009), Balakrishnan and Hashorva (2010)).
  A monograph treatment rich in applications and references is Galambos and Simonelli (2004).

When $R$ is not directly observable, but the distribution function (df) of $S$ is known, and $W$ is observable, a natural question arising from \eqref{eq:ab} is the recovery of the distribution function  of $R$, or its estimation. Such a question arises for instance while  estimating the true claim cost of a glass  insurance coverage. Indeed, if $R_i, i\ge 1$ models the losses payed to claims reported from some  glass coverage of a particular motor portfolio, the insurer is interested in the estimation of the true claim cost $W_i$. However, these costs are typically deflations of $R_i$, where the deflator $S_i$ explains the presence of fraud or other effects; in this setup $R_i$ is  not directly observable.\\
In certain cases the df of the scaling random variable is known,  or it can be estimated, which prompts the insurer to attempt to recover the df of the true losses. This is possible when the df of the random variable $W$ is a beta-product convolution, i.e., the scaling random variable
$S$ is beta distributed with positive parameters $\alpha,\beta$, see \eqref{eq:mainth:2} below. An interesting fact connected with beta-product convolutions is the characterisation of $k$-monotone functions, see Pakes and Navarro (2007), Balabdaoui and Wellner (2010).

The principal aim of this paper is the investigation of the scaling model \eqref{eq:ab} from both distributional and asymptotical interest
extending some previous findings of Hashorva and Pakes (2010) and Hashorva et al.\ (2010).
In the framework of beta random scaling we consider the inverse problems of derivation of the distribution (or the density) function of $R$ when that of $W$ is known, which can be solved by resorting to properties of the Weyl fractional-order  integral operator. In the second part of the paper we deal with the lower tail asymptotic behaviour of $W$ and $R$, dropping specific
 distributional assumptions on $S$. In particular, we investigate the min-domain of attraction of the distribution function of $W$ if that of $R$ belongs to the min-domain of attraction of some univariate extreme value distribution function.
For the beta-product convolution model we are able to derive further some converse asymptotic results.

 Distributional properties and results such as lower (upper) tail asymptotics of beta-product convolutions are of certain importance for insurance application when dealing for instance with the modeling of small (large) claims  which are typically affected by some random inflation (deflation) factor.
  In fact, from the financial point of view, insurance companies do not suffer from small claims but from the large ones.
However, understanding small claims is important for at least two reasons:
  a) claim handling is expansive even for zero-losses  or very small ones, b) the choice of deductibles and the calculation of pure premiums can be significantly improved if the effect of inflation/deflation on small claims is adequately modeled.
In finance, modeling of the effect of a deflator, which can practically ruin an investment, is very important.

We present in this paper two applications: first we investigate the asymptotic behaviour of the componentwise minima of the absolute value of elliptical random vectors where we show that the it is attracted by some multivariate df with independent components, provided that the associated random radius a regularly varying (at 0) df.
In the second application we consider the aggregation of two risks with polar representation (similar to that of
elliptically distributed risks).  Aggregation of risks is an important topic for insurance and finance,
see the recent contributions  
Dhaene et al.\ (2008), Asmussen  and Rojas-Nandaypa (2008), Embrechts  and Puccetti  (2008), Albrecher and Kortschak (2009), Geluk and Tang (2009), Mitra and Resnick  (2009), Valdez et al.\ (2009), and Dengen et al.\ (2010).

Outline of the rest of the paper: Preliminary results will be followed by Section 3 where we
discuss focus on the main distributional properties underlying the beta random scaling model.\\
Lower tail asymptotics for $W$ and $R$ related by \eqref{eq:ab} is investigated in Section 4. Two applications in Section 5 proceed the las section which contains some related results and the our proofs.

\section{Preliminaries}
In the sequel $\alpha,\beta$ are two positive constants, and $\Bab$ denotes a beta random variable
with density function
$$\frac{ \Gamma(\alpha+\beta)}{ \Gamma(\alpha)\Gamma(\beta)  }x^{\alpha-1}(1- x)^{\beta-1}, \quad x\in (0,1),$$
where  $\Gamma(\cdot)$ is the Euler gamma function.
The distribution function of $R$ will be denoted by  $H$ (abbreviate this as  $R\sim H$)
$\Fab$ is the df of $W$ with stochastic representation \eqref{eq:ab}.
with upper endpoint of $\omega\omega\in (0,\IF]$. It will be assumed that the lower endpoint of $H$ is 0 (so $H(0)=0$).
For our scaling model \eqref{eq:ab} the df of $W$ is said to be a product convolution distribution defined in terms of $H$ and the df of $S$. When $S$ is beta distributed with parameters $\alpha,\beta$ the relation between $H$ and $\Fab$, with $\Fab$ the df of $S$  is quite tractable due to the role of the Weyl fractional-order integral operator. We refer to $\Fab$ alternatively as a beta-product convolution.\\ Next, we introduce the aforementioned operator
acting on  real-valued measurable functions $h$ defined on $(0,\IF)$. For a given
constant $ \beta\in(0,\IF)$ the Weyl fractional-order integral operator  $I_\beta$ is defined  by
\BQNY (I_\beta h) (x)&=&\frac{1}{\Gamma(\beta)}\int_x^\infty
(y-x)^{\beta-1} h(y) \, dy, \quad x>0.
\EQNY
 Now, if for any $\ve>0$ we have
\BQNY
\int_\ve^\IF x^{\beta-1} \abs{h(x)}\, dx < \IF,
\EQNY
which is  abbreviated by $h\in \kal{I}_\beta$, then $(I_\beta h)(x)$ is almost surely finite for all $x\in (0,\IF)$.\\
It follows easily that
\BQN \label{PKS:14}
\HAB(x)&=&\frac{\Gamma(\alpha+ \beta)}{\Gamma(\alpha)} x^{\alpha}
(I_\beta p_{- \alpha -\beta} H)(x),\quad x\in (0,\omega), \quad \text{ with }p_{s}(x)= x^{s},s  \in (0,\IF)
\EQN
showing the importance of the Weyl fractional-order integral operator in the setup of beta random scaling.

When $\alpha=1,$ then  $\pk{B_{1,\beta}> s}=(1- s)^\beta, s\in (0,1)$. Hence \eqref{PKS:14} simplifies to
$$\overline{H}_{\alpha,\beta}(x)= \int_x^\IF(1- x/y)^\beta \, d H(y), \quad \overline{H}_{\alpha,\beta}=1- \HAB,$$
which leads us to the introduction of the Weyl-Stieltjes fractional-order integral operator $\OPJ_{\beta,g} $ with $g:(0,\IF)\to \R$
a measurable weight function   defined by
\BQNY   (\OPJ_{\beta, g} H)(x) &=&\frac{1}{\Gamma(\beta)}\int_x^{\infty}(y-x)^{\beta-1}g(y) \, d
H(y), \quad x\in (0,\omega).
\EQNY
With this notation we have
\BQN\label{eq:William}
\overline{H}_{1,\beta}(x)&=& \Gamma(\beta+1)(\OPJ_{\beta+1, p_{-\beta}} H)(x),  \quad x\in (0,\omega).
\EQN
When $\beta=d\in \N$, then $\Gamma(\beta)(\OPJ_{\beta, p_{1- \beta}} H)$ is the Williamson $d$-transform of $H$, which plays a crucial role in the analysis of $L_1$-norm Dirichlet distributions (Fang et al.\ (1990)), and Archimedean copula (McNeil and Ne$\check{s}$lehov\'{a} (2009)).

For the derivation of the lower tail asymptotics of $W$ we impose an assumption on $R$ motivated by univariate extreme value theory. Specifically, we assume that $R\sim H$ is regularly varying at 0 with some index $\gamma\in (0,\IF)$, i.e.,
\BQN \label{eq:reg}
\lim_{x\downarrow 0} \frac{ H(t x)}{ H(x)}&=& t ^ \gamma, \quad \forall t\in (0,\IF).
\EQN
Alternatively, we write  $H\in RV_\gamma$ or $R\in RV_\gamma$. Eq. \eqref{eq:reg} is equivalent with the assertion $1/R$ is regularly varying at infinity with index $-\gamma$, or the df of  $1/R$ is in the max-domain of attraction of the \FRE df $\Phi_\gamma(x)=\exp(- x^{-\gamma}), x> 0$. See Resnick (1987), Bingham et al.\ (1987),
Embrechts et al.\ (1997),  De Haan and Ferreira (2006), Jessen and Mikosch (2006), or Omey and Segers (2009)  for more details on regularly varying functions and max-domain of attractions.

\section{Distributional Properties of Beta-Product Convolutions}
In this section we discuss the stochastic model \eqref{eq:ab} with $S$ being beta distributed with parameters $\alpha$ and $\beta$. Since we assume that $H(0)=0$, then $\Fab(0)=0$, and 
$\Fab$ possesses a positive density function $h_{\alpha,\beta}$ given by
\BQN \label{PKS:14:d}
h_{\alpha,\beta} (x)&=&
\frac{\Gamma(\alpha+ \beta)}{\Gamma(\alpha) } x^{\alpha-1}  (\OPJ_{\beta, p_{-\alpha- \beta +1}} H)(x), \quad x\in (0,\omega)
\EQN
implying for any $k\inn$ 
\BQN\label{kMonoton}
h_{1,k}(x)&=& 
k \int_x^\IF (y-x)^{k-1} y^k\, d H(x) , \quad x\in (0,\omega).
\EQN
The density function  $h_{1,k}$ is a $k$-monotone function, see Balabdaoui and Wellner (2010) for recent deep results concern estimation of $k$-monotone functions. Conversely, any integrable $k$-monotone function has representation \eqref{kMonoton}, see Lemma 1 in the aforementioned paper. 

Clearly, if $H$ possesses a density function $h$, then
\BQNY
(\OPJ_{\beta, p_{-\alpha- \beta +1}} H)(x)=(I_{\beta, p_{-\alpha- \beta +1}} h)(x), \quad x\in (0,\omega).
\EQNY
In particular we have for some $\delta\in [0,1), n\inn$
\BQN \label{PKS:14:e}
h_{1,n- \delta} (x)&=&
\Gamma(n+1- \delta) (I_{n-\delta, p_{\delta- n}} h)(x), \quad x\in (0,\omega).
\EQN
Let $D^{(n)}$ denote the $n$-fold derivative operator (we write alternatively $f^{(n)}$ instead of $D^{(n)}f$ for some differentiable function $f$). If $h^{(n)}_{1,n- \delta} $ exist almost everywhere, utilising \eqref{PKS:14:e} we can recover $h$ for $\delta\in (0,1)$ as
\BQN \label{PKS:14:f}
x^{\delta - n} h(x) &=& \frac{(-1)^n}{\Gamma(n+1- \delta)} (I_\delta h^{(n)}_{1,n- \delta}) (x), \quad x\in (0,\omega),
\EQN
which follows by the properties of the Weyl fractional-order integral, see \nelem{lem:app} and Corollary 2.1 of Pakes and Navarro (2007). When $\delta=0$  \nelem{lem:app} yields further
\BQN \label{PKS:14:fb}
 h(x) &=& \frac{(-x)^n}{\Gamma(n+1)} h^{(n)}_{1,n} (x), \quad x\in (0,\omega)
\EQN
which follows also from Lemma 1 in Balabdaoui and Wellner (2010). A more general result is stated in Theorem 2.1 of the Pakes and Navarro (2007). Namely,  if $h$ exists and $h\in {\cal{I}}_{1+\alpha- \delta}$ ($h\in \cal{I}_{\alpha- \delta}$ is instead assumed therein, which is a misprint), then
\BQN \label{PKS:14:ff}
h(x) &=& (-1)^n \frac{\Gamma(\alpha)}{\Gamma(\alpha+ n- \delta)} x^{\alpha+ n- \delta -1} (I_\delta D^{(n)}(p_{1- \alpha}h_{\alpha,n- \delta})) (x), \quad x\in (0,\omega),
\EQN
provided that $h^{(n)}_{\alpha,n-\delta}$ exists almost everywhere. When $\alpha \in [0, \delta]$ formalising we arrive at:\\

\BT\label{theo:2009} Let $H, \Fab$ be as above where $\alpha\in [0, \delta], \beta=n-\delta,\delta\in [0,1)$ with  $n\in \N$.  Let $h,\hab$ denote the corresponding density functions of $H$ and $\Fab$, respectively. If  $h^{(n)}_{\alpha, \beta}$ exists almost everywhere, then \eqref{PKS:14:ff} holds. \ET


{\bf Example 1.} a) Consider the case $\alpha\in (0,\IF)$ and $\beta= d\in \N.$ If $h^{(d)}_{\alpha,d}$ exists almost everywhere and further $h\in \kal{I}_{1+ \alpha}$,  then \eqref{PKS:14:ff} implies
\BQN \label{PKS:14:ffa}
h(x) &=&  (-1)^d  \frac{\Gamma(\alpha)}{\Gamma(\alpha+d ) } x^{\alpha+d -1}  D^{(d)}(p_{1- \alpha}h_{\alpha,d}) (x), \quad x\in (0,\omega).
\EQN
b) Suppose that $\alpha=1/2$ and $\beta= d-1/2, d\in \N.$ If $h^{(d)}_{1/2,d-1/2}$ exists almost everywhere, then by \netheo{theo:2009}
\BQN \label{PKS:14:ffb}
h(x) &=&  (-1)^d  \frac{\Gamma(1/2)}{\Gamma(d) } x^{d-1}  ( I_{1/2} D^{(d)} (p_{1/2}h_{1/2,d-1/2})) (x), \quad x\in (0,\omega),
\EQN
which reduces for $\beta=1/2$ to
\BQN \label{PKS:14:ffc}
h(x) &=&  -\Gamma(1/2)  ( I_{1/2} D^{(1)}(p_{1/2}h_{1/2,1/2})) (x), \quad x\in (0,\omega).
\EQN

{\bf Example 2.} Let $\HAB$ be the df of $\Gamma_{\alpha+ \beta, \lambda}$, a Gamma random variable with positive parameters $\alpha, \lambda$ and density function given by
$$ \hab(x)=
\frac{\lambda^{\alpha}}{\Gamma(\alpha)}x^{\alpha-1}
\exp(- \lambda x), \quad x\in (0,\IF).$$
Equation \eqref{PKS:14:ff} implies that $h$ is the density function  of $\Gamma_{\alpha+ \beta, \lambda}$, a Gamma random variable with parameters $ \alpha+ \beta, \lambda $. If $\Gamma_{\alpha+ \beta, \lambda}$ is independent of $\Bab$  this means
$$ \Gamma_{\alpha, \lambda} \equaldis \Gamma_{\alpha+ \beta, \lambda} \Bab,$$
which is a well-known property of gamma and beta random variables, see e.g., Galambos and Simonelli (2004).

A key fact when dealing with independent beta products is that if $B_{\lambda, \gamma}$ is beta distributed with positive parameters $\lambda=\alpha+ \beta,\gamma$ being further independent of $\Bab$, then we have the stochastic representation
\BQN
\label{eq:leo}
  \Bab B_{ \lambda, \gamma} \equaldis B_{\alpha, \beta+ \gamma}.
  \EQN
The above stochastic representation is crucial for the recursive calculation of $\hab$. Since $h$ need not always exist, it is of some importance to recover the df $H$ when $\Fab$  is known. Utilising \eqref{eq:leo} this can be achieved iteratively as shown in our next result.
\COM{in Balakrishnan and Hashorva (2010), namely
for given $\alpha,\beta$ positive, it is possible to construct distribution functions
$H_0=H, H_1 \ldot H_{k+1}=H_{\alpha,\beta}$ determined iteratively by
\BQN\label{eq:mainth:2}
 H_{i-1}(x) &=& \frac{\Gamma( \alpha + \beta_{i} )}{\Gamma(\alpha + \beta_{i-1})}
 x^{\alpha + \beta_{i-1}} \Bigl[ (\alpha+ \beta_{i}) (I_{\delta_{i}} p_{- \alpha - \beta_{i}-1} H_{i})(x)
      -(\OPJ_{\delta_{i},  p_{- \alpha - \beta_{i}}} H_{i} )(x)\Bigr], \quad  x\in (0,\omega),
 \EQN
provided that $\beta_i,i=0\ldot k+1$ are given constants with $\beta_0=\beta, \beta_{k+1}=\beta$ and
$\delta_i= 1+ \beta_i- \beta_{i-1}\in (0,1), i=1 \ldot k+1$
}

\BT \label{eq:mainth}
Let $H,\HAB$ be two distribution functions of the random scaling model \eqref{eq:ab}. If $H(0)=0$ and $\beta_0:=\beta> \beta_1 >\cdots >
\beta_k> \beta_{k+1}:=0, k\in \{0,\N\}$ are constants such that
$\beta_{i-1}- \beta_{i} \in (0,1), i=1 \ldot k+1$, then
there exist distribution functions $H_0:=H, H_1 \ldot H_{k+1}=H_{\alpha,\beta}$ determined iteratively by
\BQN\label{eq:mainth:2}
 H_{i-1}(x) &=& \frac{\Gamma( \alpha + \beta_{i} )}{\Gamma(\alpha + \beta_{i-1})}
 x^{\alpha + \beta_{i-1}} \Bigl[ (\alpha+ \beta_{i}) (I_{\delta_{i}} p_{- \alpha - \beta_{i}-1} H_{i})(x)
      -(\OPJ_{\delta_{i},  p_{- \alpha - \beta_{i}}} H_{i} )(x)\Bigr], \quad  x\in (0,\omega),
 \EQN
with $\delta_i:= 1+ \beta_i- \beta_{i-1}$. Furthermore, $H_i, i=1 \ldot k+1$ possesses a density function $h_i$.
\ET

We illustrate next \eqref{eq:mainth:2} by two examples.

{\bf Example 3}. Consider $H, \HAB$ with $\alpha\in (0,\IF)$ and $\beta=d \in \N.$ With $\beta_i= d-i, i=0 \ldot d$ there exists $H_i$ with  density function $h_i, i=1 \ldot d$ such that $H_0=H, H_{d}=\HAB$ and
\BQN \label{eq:2010}
 H_{i-1}(x) 
&=& \frac{1}{\alpha + \beta_{i}}
 x^{\alpha + \beta_{i-1}}\Bigl[ (\alpha+ \beta_i)x^{- \alpha - \beta_{i-1}} H_{i}(x)
      - x^{- \alpha - \beta_{i}} h_{i} (x)\Bigr] \notag\\
&=& H_{i}(x)      - \frac{x h_{i}(x)}{\alpha + \beta_{i}}, \quad  x\in (0,\omega), \quad i\in \{1 \ldot d\}.
 \EQN
Consequently, if  $h^{(d)}_{\alpha,d}$ exists, then we can calculate $h$ recursively by
 $$ h_{i-1}(x)= \frac{\alpha + \beta_{i}- 1}{\alpha + \beta_{i}}h_i(x)- \frac{   xh_{i}^{(1)}(x) }{\alpha + \beta_{i}}, \quad x\in (0,\IF),$$
 which is an alternative calculation to \eqref{PKS:14:ffa}.
 Note that if $\alpha=1,$ then $H$ can be determined explicitly by the inverse of the Williamson  $d$-transform, see Proposition 3.1 in McNeil and Ne$\check{s}$lehov\'{a} (2009), or  Lemma 1 in Balabdaoui and Wellner (2010).

{\bf Example 4}.  In a financial context assume that an investment $R\sim H$ (positive) is being subjected to some deflation effect such  that the
return after a period of time (say a year) is $W\equaldis R S$  with deflator $S$ being uniformly distributed on the interval  $(0,1)$. The fact that
$S \equaldis B_{1,1}$ and \eqref{eq:mainth:2} 
imply that $H$ and the df $H_{1,1}$ of $W$ are related by
\BQN\label{bI}
 H(x) &=& H_{1,1}(x) - x h_{1,1}(x), \quad \forall x\in (0,\omega),
 \EQN
with $\omega\in (0,\IF]$ the upper endpoint of $H$. Furthermore, if $h_{1,1}^{(1)}$ exists, then almost surely in $(0,\omega)$
\BQNY
 h(x) &=& x h_{1,1}^{(1)}(x) .
 \EQNY
Consequently, we have
$$ \lim_{ x\downarrow 0} x h_{1,1}(x)=\lim_{ x \uparrow  \omega} x h_{1,1}(x)=0.$$
Next, if for some constant $\gamma$
$$\lim_{x \downarrow  0} \frac{x h_{1,1}(x)}{H_{1,1}(x)}= \gamma\in [0,1],$$
then by Proposition 2.5 in Resnick (2007) $H_{1,1} \in RV_\gamma$. Further, \eqref{bI} implies
$$ \limit{x}  \frac{ H(1/x)}{H_{1,1}(1/x)} = 1- \gamma,$$
which can be also written alternatively as (set $R^*:=1/R$)
$$ \limit{x}  \frac{\pk{R^* \xi> x }}{ \pk{R^*> x}} =\frac{1}{ 1- \gamma},$$
where $\xi= 1/ B_{1,1}$ is a Pareto random variable with parameter 1. Karamata's Theorem (see e.g., De Haan and Ferreira (2006), Resnick (2007))
yields   thus  if $\gamma\in [0,1)$, that  also $H\in RV_\gamma$.
 Another  proof of this fact is given in Proposition 5.2 of Maulik and Resnick (2004). By \eqref{bI} a converse result can be easily established. Note in passing that since $h(x)=x h_{1,1}^{(1)}(x)$ we have $h\in RV_{\gamma- 1}$ if and only if $h_{1,1}^{(1)} \in RV_{\gamma- 2}, \gamma\inr$.

\section{Lower Tail Asymptotics}
In the recent contributions Hashorva and Pakes (2010), Hashorva et al.\ (2010) discuss the asymptotic behaviour of the survival function $\overline{H}_{\alpha,\beta}$
assuming that $H$ belongs to some max-domain of attraction of a univariate extreme value df. \\
Hashorva and Pakes (2010) shows that $H$ and $\HAB$ belong (if so) to the same max-domain of attraction. The practical importance of these findings is that for insurance and finance random scaling models do not allow the deflator to change the max-domain of attractions of the random payment.
Interesting asymptotic results for our random scaling models can be found in the context of Archimedean copula
in  Charpentier and Segers (2007, 2008, 2009), see Remark \eqref{rem:th1} below. For such copulas in the aforementioned papers the asymptotics of their density generator $\psi$ at 0 is derived. The connection with our scaling model is immediate since $\psi$ equals the survival function  $\overline{H}_{1,\beta}, \beta\inn$.

Complementing the findings of Hashorva and Pakes (2010) we focus next on the lower tail asymptotics of $\HAB$, which boils down to determination of the min-domain of attractions of  $\HAB$.\\
When dealing with positive random variables the  min-domain of attraction, for say the df $H$, is determined by the max-domain of attraction of
the df $H_*$ of $1/R$.
 Since $H$ is a df with lower endpoint 0,
only the \FRE or the Gumbel max-domain of attraction for $H_*$ is possible. The first assumption (\FRE) is equivalent with
$H$ satisfying \eqref{eq:reg} with some positive index $\gamma$. Actually, this situation is simple
since by Lemma 4.1 and Lemma 4.2 in Jessen and Mikosch (2006) 
 if $H$ satisfies  \eqref{eq:reg} with some $\gamma\in (0,\IF)$, and $S\in RV_\alpha$ with $\alpha\in (0,\IF)$,
 then
 \BQN \label{eq:neu}
\Fab\in RV_{\gamma^*}, \quad \gamma^*:= \min(\gamma, \alpha).
 \EQN
 Note in passing that when  $\alpha=\gamma$ the asymptotic behaviour of random product follows from the well-known result of Embrechts and Goldie (1980). \\
 The beta-product convolution model is included in the above assumption since if  $S\sim \Bab$ then
\BQN \label{tail:beta}
\pk{ S< s}
&=& (1+o(1))\frac{\Gamma(\alpha+ \beta)}{\Gamma(\alpha+1)\Gamma(\beta)} s^{\alpha}, \quad s\downarrow 0
\EQN
implying $S\in RV_\alpha$.

In the next theorem we consider initially the case $H_*$ is in the Gumbel max-domain of attraction, and then prove a converse asymptotic
result for the \FRE case.\\

\BT \label{th:1}Let $R\sim H$ be a positive random variable which is independent of $S=\Bab,$ and define the beta-product convolution df $\Fab$ with lower endpoint 0 via the random scaling model \eqref{eq:ab}. We have:

\COM{Furthermore, if $\alpha\not=\gamma$
\BQN\label{eq:ref:in}
\lim_{ s \downarrow 0} \frac{  \Fab(s)}{ H(s)}&=& K_{\alpha,\gamma} \in (0,\IF), \quad
\EQN
where $K_{\alpha,\gamma} =\E{R^{-\alpha}}$ if $\alpha>\gamma$ and
$$ K_{\alpha,\gamma}= \frac{\Gamma(\alpha+ \beta)\Gamma(\alpha- \gamma )}{\Gamma(\alpha)\Gamma(\alpha+ \beta- \gamma)}, \quad \alpha\in (0,\gamma).$$
}

b) If $1/R$ has df in the Gumbel max-domain of attraction, then $\Fab\in RV_{\alpha}$.

c) If $\Fab$ satisfies \eqref{eq:reg} with some $\gamma\in (0,\alpha)$, then $H\in RV_\gamma$ and $h_{\alpha,\beta}\in RV_{\gamma- 1}$.\\
\ET

\begin{remark} \label{rem:th1}
\end{remark}
\COM{
1. If $\gamma$ in \netheo{th:1} equals $\alpha,$ and for some $c\in (0,\IF)$
$$ \pk{1/R> u}= (1+o(1)) c x^{- \gamma} \, \quad u\to \IF,$$
then  \eqref{eq:ref:in} holds with $K_{\alpha,\gamma}:= $. This follows by Exercise in Resnick (2007).
}
1)  In the setup of Archimedean copula Charpentier and Segers (2009) consider the asymptotics at 0 of $\psi^{-1}$ with $\psi$ the generator of some  Archimedean copula (in the notation of McNeil and Ne$\check{s}$lehov\'{a} (2009)). In the light of the findings of
the aforementioned paper $\psi^{-1}= (\overline{H}_{1,1})^{-1}$, i.e., it is the inverse of the survival function of a beta-product convolution ($\alpha=1,\beta=1$). By Proposition 2.6 (v) in Resnick (2007) $\psi^{-1} \in RV_{\alpha}, \alpha\in [0,\IF]$ implies that $\overline{H}_{1,1}$ is
regularly varying at infinity with index $1/\alpha$. For the general $k$-dimensional Archimedean copula  $\psi^{-1}= (\overline{H}_{1,k-1})^{-1}, k\in \N$. Consequently, the findings of Charpentier and Segers (2009) concern the asymptotics of the survival functions $\overline{H},\overline{H}_{1,k}$ and that of $h_{1,1}(x)$ as $x\to \IF$. Note further that
the identity \eqref{bI} of Example 4 can also be utilised to deal with these functions.

\COM{ Note that when $\alpha=\gamma,$ and $\E{S^{-\alpha}}=\E{R^{-\alpha}}=\IF$, then again \eqref{eq:ref:in}  holds and further
$K_{\alpha,\alpha}=\IF$, see Lemma 4.1 in Jessen and Mikosch (2006), and Problem 7.8 in Resnick (2007). In the special case that
$$F(x)= (1+o(1))c x^{\alpha}, \quad G(x)= (1+o(1))c x^{\alpha},\quad x \downarrow 0$$
for $c,\alpha \in (0,\IF)$ we have
$$ K_{\alpha,\alpha}= - (1+o(1))\alpha c^{2 \alpha} x^\alpha \ln x, \quad x \downarrow 0,$$
which follows from the aforementioned lemma.
}

2) If $\alpha\in (0,\IF), \beta=d$, then by \eqref{PKS:14:ffa} regular variation of $h$ holds if the same is true for $D^{(d)}(p_{1- \alpha} h_{\alpha,d})$.
Clearly, for $\alpha=1$ the latter reduces to $h_{\alpha,d}^{(d)}$. %

3) For any $\alpha,\lambda \in (0,\IF)$ and $S$ a positive scaling random variable $S^\lambda \in RV_{\alpha}$ is equivalent with
$S \in RV_{\alpha\lambda}$. Consequently, our asymptotic results above apply also when $S^\lambda$ is a beta random variable.

\section{Applications}
In this section we provide two applications. Motivated by the findings of Kabluchko (2010) we derive first the joint asymptotic independence of sample minima considering bivariate elliptical random vectors, which was shown
in the aforementioned paper for the special case of Gaussian random vectors.

Our second application is concerned with lower tail asymptotics for polar distributions which is also related to the lower tail asymptotics of aggregated risks. As mentioned in the Introduction aggregation is a central topic in various applications; for insurance and financial applications see e.g., Denuit et al.\ (2005) and Dengen et al.\ (2010). \\

\subsection{Asymptotics of Minima for Elliptical Samples}
When $\vk{U}$ is uniformly distributed on the unit sphere of $\R^k, k\ge 2$, and $A$ is a $k$-dimensional nonsingular real matrix, then the random vector $\X \equaldis  R A \vk{U}$ is elliptically distributed. It is well-known (cf.\ Fang et al.\ (1990)) that the distribution function of $\X$ depends on $\SI:=AA^\top$ but not on the matrix
$A$ itself. In view of the properties of $\vk{U}$ (cf.\ Fang et al.\ (1990)), if further
the main diagonal of $\Sigma$ consists of 1s, i.e., $\Sigma$ is a correlation matrix, then by Lemma 6.1 in Berman (1983)
 \BQN \label{eq:el2}
 X_i\equaldis X_1\equaldis - X_1\equaldis R U_1, \quad 1  \le i \le k,
 \EQN
with $U_1$ the first component of $\vk{U}$.  Furthermore, we have the stochastic representation
 $$X_1^2 \equaldis R^2 B_{1/2, (k-1)/2}, $$
 where $R$ is independent of $B_{1/2, (k-1)/2}$. Clearly, the random variable $\abs{X_1}$ is a deflation of $R$ by  $S=\sqrt{B_{1/2, (k-1)/2}}$.\\
Next, \eqref{eq:neu} and the above stochastic representation (recall also \eqref{tail:beta}) imply
 $\abs{X_{11}}$ has df  $Q\in RV_\gamma, \gamma \in (0,1]$
if for instance $H\in RV_\gamma, \gamma \in (0,1]$.  For such df $Q$ we define constants
$a_n,n\ge 1$ asymptotically by
\BQNY 2n\pk{a_n^{-1} \ge X_{11} > 0}&=& 1.
\EQNY
It follows that $a_n^{-1}=L(1/n) n^{-\gamma},$ with $L\in RV_0$ a positive slowly varying function at 0, i.e., $\limit{n} L(c/n)/L(1/n)=1, \forall c\in (0,\IF)$.
For such constants we have the convergence in distribution as $n\to \IF$
\BQNY
a _n M_{ni} \todis  \kal{M}_i \sim \kal{G}_\gamma, \quad M_{ni}:= \min_{1 \le j \le n} \abs{X_{ji}}, \quad 1 \le i\le k,
\EQNY
with  $\kal{G}_\gamma$ given by
\BQN\label{GY}
\kal{G}_\gamma(x)&=& 1- \exp(-x^\gamma), \quad x> 0.
\EQN
Since for  $u\in (0,\IF)$ small enough
$$\pk{\abs{X_{1i}}< u, \abs{X_{1j}}< u}=0, \quad 1 \le i \not=j \le k$$
we have
$$\lim_{u\downarrow 0} \frac{ \pk{\abs{X_{1i}}< u, \abs{X_{1j}}< u}}{\pk{\abs{X_{11}}\le u}}=0, \quad 1 \le i \not=j \le k.$$
Consequently, the sample minima $\vk{M}_n=(M_{n1} \ldot M_{nk})$ has asymptotic independent components, meaning that the joint convergence in distribution
\BQNY
\Bigl( a _n M_{n1} \ldot  a _n M_{nk} \Bigr)  &\todis & \Bigl(\kal{M}_1  \ldot \kal{M}_k\Bigr), \quad n\to \IF
\EQNY
holds with $\kal{M}_1 \ldot \kal{M}_k$ independent with df $\kal{G}_\gamma$. \\
\COM{
\BQNY
\pk{ M_{ni} \le x/a_n}&=& 1- \pk{\abs{X_{11}} > x/a_n}^n\\
&=& 1- \exp( - \limit{n} n \pk{\abs{X_{11}}\le  x/a_n})\\
&=& 1- \exp( - \limit{n} 2 n \pk{0 < X_{11} \le x/a_n})\\
&=& 1- \exp( - 2  x^\gamma).
\EQNY
}
To this end, we note that the joint convergence in distribution above can be reformulated for the more general class of asymptotically elliptical random vectors, see Hashorva (2005) for asymptotic properties.

\subsection{Aggregation of Two Risks}
If $\X$ is a $k$-dimensional elliptical random vector as above, then for given constants $\mu_i, 1 \le i\le k$ we have
$$ \sum_{1 \le j \le k} \mu_j X_j \equaldis \sqrt{ \sum_{1 \le j \le k} \mu_j^2} X_1,$$
which is a well-known property for the Gaussian random vectors. Moreover, if $Be(q), q \in (0,1]$ denotes the df of a Bernuli random variable assuming values $-1,1$,  for any pair $X_i, X_j, i\not=j$ we have
\BQN
\label{eq:doub}
 (X_i, X_j) \equaldis \Bigl(O_1, \rho_{ij} O_1+ \sqrt{1 -\rho_{ij}^2}  O_2 \Bigr),  \quad (O_1,O_2) \equaldis (R_* T_1 S, R_* T_2
 \sqrt{1- S^2}), \quad S \sim \sqrt{B_{1/2,1/2}},
 \EQN
with $\rho_{ij} \in (-1,1)$ the $ij$th entry of the correlation matrix $\SI$, and $(O_1,O_2)$ a bivariate spherical random vector with
positive associated random radius $R_*$ such that $R_*^2 \eqd R^2 B_{1, k/2- 1}$ and $T_i \equaldis Be(1/2), i=1,2$.
 Furthermore $T_1,T_2, R_*, S$ are mutually independent.

Since we are concerned with asymptotic results, the distribution assumption on $S$ above can be dropped.
We consider next the lower tail asymptotics of a bivariate polar random vector $(X,Y_\rho), \rho\in (-1,1)$ with stochastic representation
\BQN
\label{eq:doub:2}
 (X, Y_\rho) \equaldis \Bigl(T_1 R S, \rho T_1 R S+ \tilde \rho T_2 R \sqrt{1- S^2}\Bigr), \quad \tRO:=\sqrt{1- \rho^2},
 \EQN
where $T_i\eqd Be(q_i), q_i\in (0,1]$, $R\sim H$ and $S\sim G$ such that
$$G(0)=H(0)=0, \quad G(1)=1.$$
 As in the elliptical setup here again $T_1,T_2, R, S$ are assumed to be mutually independent. Since $\abs{X}\equaldis R S,$ then the lower tail asymptotics of $\abs{X}$ can be established by \netheo{th:1} under asymptotic assumptions on both $R$ and $S$.\\
We note in passing that asymptotic properties of $(X,Y_{\rho})$ with $\rho\in (-1,1)$ random are discussed in
the recent contribution Manner  and Segers (2009). Further, remark that as mentioned in Remark \ref{rem:th1}, we do not need to specify the df $G$, apart from  the asymptotic condition  $G\in RV_\alpha$.\\
We derive next the lower tail asymptotics of $\abs{Y_\rho}$ under the following additional assumption: For all positive $t$ small enough
\BQN\label{eq:G}
G(\rho + t)- G(\rho-t)=L_\rho(t) t^{\garho},\quad G(\tRO + t)- G(\tRO-t)=L_{\tRO}(t) t^{\gatrho}, \quad \garho,\gatrho\in [0,\IF),
\EQN
with $L_\rho,L_{\tRO}\in RV_0$ two positive functions. Condition \eqref{eq:doub:2} can be easily checked. In the special case that $G$ possesses a  positive density function  $g$ continuous at $\rho$  and $\tRO$ condition \eqref{eq:G} is satisfied with
\BQN \label{eadis:c}
\garho=\gatrho=1, \quad \text{ and } L_\rho(t)= (2+o(1)) g(\rho), \quad L_{\tRO}(t)= (2+o(1)) g(\tRO), \quad t\downarrow 0.
\EQN
We have now the following result.

\BS \label{last} Let $(X,Y_\rho), \rho\in (0,1)$ be a bivariate random vector with stochastic representation \eqref{eq:doub:2}. Suppose
that
$$R\sim H\in RV_{\gamma},  \quad S\sim G\in RV_{\alpha}, \quad \alpha,\gamma \in (0,\IF)$$
 and \eqref{eq:G} holds with some $\garho, \gatrho$ and $L_\rho, L_{\tRO}$. Assume further that when $\alpha_{\rho}=\gatrho,$ then
$L_\rho(x)=cL_{\tRO}(x), \forall x>0$ with some positive constant $c$. Then
 \BQN
 \abs{X} \in RV_{\gamma_1}, \quad \abs{Y_\rho} \in RV_{\gamma_2},
 \EQN
 where $\gamma_1=\min(\alpha,\gamma) $ and $ \gamma_2=\min(\gamma, \alpha_{\rho},\alpha_{\tRO}).$
\ES
A simple instance for which we can apply \neprop{last} is when $G$ possesses a continuous positive density function $g$. In view of
\eqref{eadis:c} the index $\gamma_2$  equals $\min(\gamma,1)$. The following corollary is an immediate consequence of the above discussion.\\

\BK\label{korr:l} Let $G,H, (X,Y_\rho), \rho\in (-1,1)$ be as in  \neprop{last}, and let $(X_n,Y_n), n\ge 1$ be independent  bivariate random vectors with
the same df as $(X,Y_\rho)$. If $G\in RV_{1}$, and it possesses a continuous positive density function $g$, then we have the joint convergence in distribution
\BQN
\Bigl( a_n \min_{1 \le j \le n} \abs{X_j},  b_n \min_{1 \le j \le n} \abs{Y_j}\Bigr) & \todis &  \Bigl( \kal{M}_1, \kal{M}_2\Bigr) \quad n\to \IF,
\EQN
with $\kal{M}_1, \kal{M}_2$ independent with df  $\kal{G}_{ \min(\gamma,1)}$ defined in \eqref{GY} and constants
 $a_n,b_n,n\ge 1$ satisfying $\pk{\abs{X} < 1/a_n}= \pk{\abs{Y_\rho}< 1/b_n}=1/n$ for all large $n$.
\EK

\def\SRO{S_\rho}

\section{Further Results and Proofs}
Next we present first two lemmas and then proceed with the proofs of the claims in the previous sections.

\BL \label{lem:app}
Let $\beta,c$ be positive constants, and let $h:(0,\IF) \to R $ be a given positive measurable function. \\
a) If for some $r\in (0,\IF)$ we have $\int_0^\infty x^{r-1} \abs{h(x)}\, dx < \IF $, then $I_\beta h$ is continuous at 0, i.e.,
\BQN\label{eq:ap:I0}
\lim_{\beta \downarrow 0} I_\beta  h &=& I_0 h =h.
\EQN
b) If $h\in \kal{I}_{\beta}$ and when $\beta\in (0,1)$
\BQNY
\int_x^{\delta} (y-x)^{\beta -1} h(y) \, dy < \IF, \quad x> 0
\EQNY
holds for some $\delta> x$, then $(I_\beta h)(x)$ is finite and continuous for all $x>0$.

c) If $h\in \kal{I}_{\beta+ c}$, then
\BQN
I_\beta  I_c h =   I_c  I_\beta h = I_{\beta +c} h.
\EQN
d) If the $n$-fold derivative $D^{(n)}h $ exists
almost everywhere and $D^{(n)}h  \in \kal{I}_\beta$, then
\BQN\label{DN}
   D^{(n)} I_\beta h &=& I_\beta D^{(n)} h.
   \EQN
e) For any df  $H$ with $H(0)=0$  and upper endpoint $\omega\in (0,\IF]$ we have
\BQN\label{eeDN}
(\OPJ_{\beta+1, p_{1- \beta}} H)(x)&=& x (I_{\beta}p_{-1 - \beta} \overline{H})(x), \quad x\in (0,\omega).
\EQN
\EL
\prooflem{lem:app} By the assumption on $h$  its Mellin transform exists, consequently
utilising statement (d) Lemma 2.2 in Pakes and Navarro (2007) establishes the proof.
Statement b) is mentioned in the Introduction of Hashorva and Pakes (2010). Both $c)$ and $d)$ are shown in Lemma 8.1 of the aforementioned paper.
\eqref{eeDN} follows immediately from \eqref{PKS:14} and \eqref{eq:William}, and thus the proof is complete. \QED

\begin{lem} \label{lem:11}
Let $T_1,T_2$ be two random variables taking values $-1,1$ with
$\pk{T_1T_2=-1}\in (0,1]$ being independent of the scaling random variable $S\sim G$ with $G(0)=1- G(1)=0$. For given $\rho\in (0,1)$ set $\SRO:= \abs{\rho T_1 S+ \tRO T_2  \sqrt{1- S^2}}$. If $G$ satisfies \eqref{eq:G}, then we have
\BQN\label{eq:lem:11}
\pk{\SRO \le u}&=&(1+o(1))q_{1,-1} (\rho u)^{\gatrho} L_{\tRO }(u) +
(1+o(1))q_{-1,1} (\tRO  u)^{\garho} L_{\rho }(u), \quad u\downarrow 0,
\EQN
where $q_{i,j}:=\pk{T_1=i,T_2=j}, i,j\in \{-1, 1\}$.
\end{lem}
Note in passing that  if $G$ possesses a positive density function $g$ continuous at $\rho$ and $\tRO$, then \eqref{eq:lem:11} reduces to
\BQNY
\pk{\SRO\le u}&=& (1+o(1))2 \pk{T_1 T_2=-1}[g(\rho) \tRO + g(\tRO) \rho] u, \quad u \downarrow 0.
\EQNY

\prooflem{lem:11} By the assumptions  $S\in (0,1)$ almost surely, and $T_j,j=1,2$ assumes only two values $\{-1,1\}$. Hence we may write
for any $u\in (0,1)$ small enough
\BQNY
\pk{\SRO\le u}&=&  \pk{T_1=1, T_2=-1} \pk{\abs{\rho S- \tRO \sqrt{1- S^2}}\le u}\\
&&+ \pk{T_1=-1, T_2=1} \pk{ \abs{\tRO \sqrt{1- S^2}- \rho S}\le u}.
\EQNY
Using further the fact that $S$ is independent of $T_1,T_2$ we obtain
\BQNY
\pk{\abs{\rho S- \tRO \sqrt{1- S^2}}\le u}&=&
\pk{- u\le \rho S- \tRO \sqrt{1- S^2}\le u}\\
&=& (1+o(1))\int_{\tRO - (1+o(1)) \rho u}^{\tRO + (1+o(1)) \rho u} \, d G(s)\\
&=& (1+o(1))(\rho u)^{\gatrho} L_{\tRO }(u), \quad u \downarrow 0.
\EQNY
As above we have further
\BQNY
\pk{\abs{\tRO \sqrt{1- S^2}- \rho S}\le u}&=& (1+o(1))(\tRO  u)^{\gamma_{\rho}} L_{\rho}(u), \quad u \downarrow 0,
\EQNY
thus the result follows. \QED

\prooftheo{eq:mainth} Theorem 3.3 in Hashorva and Pakes (2010) shows an iterative formula for calculating the survival function $\overline{H}$
 when the survival function  $\overline{H}_{\alpha,\beta}$ is known.
 Our proof here is established with the same arguments of the aforementioned theorem utilising further \eqref{eq:leo}. \QED

\prooftheo{th:1}  a) It is well-known (see e.g., Resnick (1987)) that if $R^*=1/R$ has df  in the Gumbel max-domain of attraction, then $\E{(R^*)^\delta}< \IF$
for any $\delta>0$. Hence since $\Bab$ is regularly varying at 0 with index $\alpha$ the claim follows by Breiman's Lemma
(see  e.g., Jessen and Mikosch (2006), Denisov  and  Zwart (2007),  or Resnick (2007)).

b) We show the proof utilising \netheo{eq:mainth}. With the notation of the aforementioned theorem
there exist distribution functions $H_0:=H, H_1 \ldot H_{k+1}=H_{\alpha,\beta}$ determined iteratively by
\BQNY
 H_{i-1}(x) &= & \frac{\Gamma( \alpha + \beta_{i} )}{\Gamma(\alpha + \beta_{i-1})}
 x^{\alpha + \beta_{i-1}} \Bigl[ (\alpha+ \beta_{i}) (I_{\delta_{i}} p_{- \alpha - \beta_{i}-1} H_{i})(x)
      -(\OPJ_{\delta_{i},  p_{- \alpha - \beta_{i}}} H_{i} )(x)\Bigr], \quad \forall x\in (0,\IF),
 \EQNY
with $\delta_i:= 1+ \beta_i- \beta_{i-1}\in [0,1)$ and $\beta_0:=\beta> \beta_1 >\cdots >
\beta_k> \beta_{k+1}:=0$. By the assumption on $\Fab$ for any $x>0$ we have  (set $\lambda_{k+1}:=\alpha +\beta_{k+1}+1)$
\BQNY
\Gamma(\delta_{k+1}) (I_{\delta_{k+1}} p_{- \lambda_{k+1}} H_{k+1})(x)&=& \int_{x}^\IF
(y- x)^{\delta_{k+1}-1} y^{- \lambda_{k+1}} H_{k+1}(y)\, dy\\
&=& x^{\alpha- \beta_{k}} H_{k+1}(x)\int_{1}^\IF(y- 1)^{\delta_{k+1}-1} y^{- \lambda_{k+1}} H_{k+1}( x y)/H_{k+1}(x)\, dy.
\EQNY
By Karamata's Theorem (see Embrechts et al.\ (1997), or Resnick (2007)) 
\BQNY
\Gamma(\delta_{k+1})(I_{\delta_{k+1}} p_{- \lambda_{k+1}} H_{k+1})(x)&=& (1+o(1))
 x^{- \alpha- \beta_{k}} H_{k+1}(x) \int_{1}^\IF(y- 1)^{\delta_{k+1}-1} y^{- \lambda_{k+1}+ \gamma}\, dy\\
  &=& (1+o(1))x^{- \alpha- \beta_{k}} H_{k+1}(x) \frac{\Gamma(\delta_{k+1})
 \Gamma(\lambda_{k+1}- \delta_{k+1}-  \gamma)}{\Gamma(\lambda_{k+1} -\gamma )}, \quad x \downarrow 0.
\EQNY
Similarly we obtain
\BQNY
\Gamma(\delta_{k+1})(\OPJ_{\delta_{k+1},  p_{ 1- \lambda_{k+1}}} H_{k+1} )(x) &=& \int_{x}^\IF
(y- x)^{\delta_{k+1}-1} y^{1- \lambda_{k+1}} \, dH_{k+1}(y)\\
&=&x^{1- \lambda_{k+1}} x^{\delta_{k+1}-1}\int_{1}^\IF(y- 1)^{\delta_{k+1}-1} y^{1- \lambda_{k+1}} \, d H_{k+1}( x y)\\
&=& (1+o(1))\gamma (I_{\delta_{k+1}} p_{- \lambda_{k+1}} H_{k+1})(x), \quad x \downarrow 0.
\EQNY
Consequently,
\BQNY
 H_{k}(x) &= &  (1+o(1))\frac{\Gamma( \alpha + \beta_{k+1} )}{\Gamma(\alpha + \beta_{k})}
 (\alpha+ \beta_{k+1} - \gamma ) \frac{
 \Gamma(\alpha+ \beta_{k}-  \gamma)}{\Gamma(\alpha + \beta_{k+1}+1 -\gamma )} H_{k+1}(x)\\
&= &  (1+o(1))
 \frac{\Gamma(\alpha+ \beta_{k+1})
 \Gamma(\alpha+\beta_{k}-  \gamma)}{\Gamma(\alpha + \beta_{k}) \Gamma(\alpha + \beta_{k+1} -\gamma )} H_{k+1}(x), \quad x \downarrow 0,
 \EQNY
hence $H_k\in RV_\gamma$. Proceeding iteratively we find that $H_0=H \in RV_\gamma$.\\
Next, in view of \eqref{PKS:14:d} if  $H\in RV_\gamma, \gamma  \in (0,\IF)$, then  as above we obtain
\BQNY \label{PKS:14:d}
h_{\alpha,\beta} (x)&=&\frac{\Gamma(\alpha+ \beta)}{\Gamma(\alpha) \Gamma(\beta)} x^{\alpha-1}
\int_x^\IF(y -x) ^{\beta -1} y ^{-\alpha- \beta +1} \ d H(y)\\
&=&\frac{\Gamma(\alpha+ \beta)}{\Gamma(\alpha) \Gamma(\beta)} H(x)
\int_1^\IF(y -1) ^{\beta -1} y ^{-\alpha- \beta +1} \ d H(x y)/H(x)\\
&=&(1+o(1))\gamma  \frac{\Gamma(\alpha+ \beta)}{\Gamma(\alpha) \Gamma(\beta)} \frac{H(x)}{x}
\int_1^\IF(y -1) ^{\beta -1} y ^{-\alpha- \beta + \gamma}\, dy,  \quad x\downarrow 0\\
&=&(1+o(1))\gamma  \frac{\Gamma(\alpha+ \beta)}{\Gamma(\alpha) }
\frac{\Gamma(\alpha - \gamma)}{\Gamma(\alpha+ \beta - \gamma)}
\frac{H(x)}{x}\\
&=& \gamma \frac{\Fab}{x}, \quad x \downarrow 0
\EQNY
 establishing thus the claim. \QED

\proofprop{last} The proof follows with the same arguments as in the proof of \netheo{th:1} applying further the result of \nelem{lem:11}. \QED


\bibliographystyle{plain}

\end{document}